\documentclass[12pt,journal, onecolumn, draftcls]{IEEEtran}

\usepackage{verbatim}%comment
\usepackage{amsmath,graphicx}
\usepackage{mathrsfs}
\usepackage{cite}%bibliography
\usepackage{amsthm}%theorenm, lemma, definition, ...
\usepackage{epsfig,psfrag}%subfigure,enumerate
\usepackage{amsfonts,amssymb}%math fonts and symbols

\begin{document}
% can use linebreaks \\ within to get better formatting as desired
\title{Cram\'{e}r-Rao Bound for Sparse Signals Fitting the
Low-Rank Model with Small Number of Parameters}

\author{Mahdi~Shaghaghi and~Sergiy~A.~Vorobyov% <-this % stops a space
\thanks{M. Shaghaghi is with the Department
of Electrical and Computer Engineering, University of Alberta,
Edmonton, AB, T6G 2V4 Canada (e-mail: mahdi.shaghaghi@ualberta.ca).
S.~A.~Vorobyov is with Aalto University, Department of Signal
Processing and
Acoustics, Finland.}% <-this % stops a space
%\thanks{This work is supported in part by the Natural Sciences and Engineering
%Research Council (NSERC) of Canada.}% <-this % stops a space
}

% The paper headers
%\markboth{IEEE TRANSACTIONS ON SIGNAL PROCESSING,~Vol.~6, No.~1, January~2007}%
%\markboth{IEEE SIGNAL PROCESSING LETTERS}%
%{Shell \MakeLowercase{\textit{et al.}}: Bare Demo of IEEEtran.cls
%for Journals}
% The only time the second header will appear is for the odd numbered pages
% after the title page when using the twoside option.
%
% *** Note that you probably will NOT want to include the author's ***
% *** name in the headers of peer review papers.                   ***
% You can use \ifCLASSOPTIONpeerreview for conditional compilation here if
% you desire.

% make the title area
\maketitle

\begin{abstract}
%\boldmath
In this , we consider signals with a low-rank covariance
matrix which reside in a low-dimensional subspace and can be written
in terms of a finite (small) number of parameters. Although such 
signals do not necessarily have a sparse representation in a finite 
basis, they possess a sparse structure which makes it possible to 
recover the signal from compressed measurements. We study the statistical 
performance bound for parameter estimation in the low-rank signal 
model from compressed measurements. Specifically, we derive the
Cram\'{e}r-Rao bound (CRB) for a generic low-rank model and we show that 
the number of compressed samples needs to be larger than the number of
sources for the existence of an unbiased estimator with finite
estimation variance. We further consider the applications to
direction-of-arrival (DOA) and spectral estimation which fit into the
low-rank signal model. We also investigate the effect of compression
on the CRB by considering numerical examples of the DOA estimation
scenario, and show how the CRB increases by increasing the 
compression or equivalently reducing the number of compressed 
samples.

\end{abstract}


% For peerreview papers, this IEEEtran command inserts a page break and
% creates the second title. It will be ignored for other modes.
\IEEEpeerreviewmaketitle

\section{Introduction}
Signals with sparse representations can be recovered from much 
less number of measurements than the number of samples given by 
the Nyquist rate using compressed sensing (CS) methods
\cite{Donoho06:CS,Candes06:ExactCS,Candes06:CS}. Such measurements
can be obtained by correlating the signal with a number of sensing
waveforms \cite{Candes08:CS,Laska07:CS,Vorobyov11:CS,Vorobyov14:CS,Xam1,Xam2}. 
The algorithms used for recovering the signals from such
measurements exploit the sparsity of the signals in a proper basis 
(see \cite{Candes06:CS,FigNowak07,Al1,Al2,Al3,Al4,Al5,Al6,Al7} to 
mention just a few existing algorithms).

There are signals which inherently possess a sparse structure
meaning that they can be defined by a small number of parameters.
However, such signals may not necessarily be represented as sparse
signals using a proper finite basis, i.e., there may not exist or 
be known a finite basis such that the transformation of the signal 
to that basis results in a small number of non-zero coefficients. 
For example, consider a signal composed of a linear combination 
of sinusoids. Such a signal generates sparse coefficients by the
discrete-time Fourier transform (DTFT), but its representation in
the Fourier basis obtained by the discrete Fourier transform (DFT)
exhibits frequency leakage \cite{MSH11:Nested_LS}. Although the DTFT
is a proper transformation, as it results in a small number of
non-zero coefficients for the considered signal, it is not a finite
basis and cannot be used in conventional CS recovery methods which
rely on a finite sparsity basis. Such methods have poor performance
for the considered signals if the DFT basis is used
\cite{Scharf11:Basis_Mismatch}.
In this , we consider a general class of sparse signals which
are represented by a small number of parameters in a low-rank signal
model. Our goal is to study the performance bounds for the
estimation of unknown parameters and also the reconstruction of this
class of signals from compressed measurements.

The Cram\'{e}r-Rao bound (CRB) \cite{VanTrees} for estimating a sparse 
parameter vector from compressed measurements has been studied in
\cite{Eldar10:Sparse_CRB}. However, the signal model in
\cite{Eldar10:Sparse_CRB} considers signals which can be represented
by a finite sparsity basis. Then, the CRB is computed using
approaches from the theory of constrained CRB in
\cite{Gorman90:Constrained_CRB} and \cite{Stoica98:Constrained_CRB}. 
The constrained CRB for estimating a low-rank matrix from compressed
measurements has been studied in \cite{Nehorai11:LowRank_CRB}. In
this , we consider a different signal model which does not
involve the constraint on the rank of a matrix. The CRB for
parameter estimation in compressed sensing has been also studied in
\cite{Scharf13:CS_Param_CRB,Ramasamy12:CS_CRB_asilomar,Ramasamy14:CS_param_TSP}. 
In \cite{Scharf13:CS_Param_CRB}, the signal of interest is assumed
to be a function of real-valued parameters, and it is not assumed to
be necessarily sparse in a finite basis. The CRB is computed and
bounded for different realizations of the measurement matrix. The
signal model considered in \cite{Ramasamy12:CS_CRB_asilomar} and
\cite{Ramasamy14:CS_param_TSP} is different from the one studied 
in this  in two aspects. Firstly, in 
\cite{Ramasamy12:CS_CRB_asilomar} and \cite{Ramasamy14:CS_param_TSP}, 
a noiseless signal is first compressed and then white noise is added
to the compressed signal. In contrast, we first add the noise to the
signal and then the result is compressed. This results in a
different distribution for the compressed measurements. Secondly and 
more significantly, in \cite{Ramasamy12:CS_CRB_asilomar} and
\cite{Ramasamy14:CS_param_TSP}, the signal is a vector which depends 
on a number of parameters, whereas in this , the signal is 
composed of a parametrized matrix multiplied by a vector of 
coefficients. This structure of the signal enables us to
derive a closed-form expression for the CRB of the parameters.

In this , we extend the results of \cite{Scharf13:CS_Param_CRB} 
for a low-rank signal model. We derive the CRB for real and 
complex-valued parameters. Furthermore, multiple signal snapshots 
are considered, whereas in \cite{Scharf13:CS_Param_CRB}, the signal 
model consists of only a single signal snapshot. We also study the 
minimum number of compressed samples required for unbiased 
estimation with finite variance. Furthermore, the applications to 
direction-of-arrival (DOA) and spectral estimation which fit into 
the low-rank signal model are also studied. Finally, numerical 
examples for the DOA estimation problem are given to illustrate 
the effect of compression on the CRB.

\section{Signal Model}
\label{sec:Sys_Mod}
Consider the signal $\boldsymbol{x}(t) \in \mathbb{C}^{N_x\times 1}$
at time instant $t$ to be of the form
\begin{equation}
\boldsymbol{x}(t) = \boldsymbol{A} \boldsymbol{d}(t)
\label{eq:matrix_sig_mod}
\end{equation}
where $\boldsymbol{A} \in \mathbb{C}^{N_x\times K}$ is a tall matrix
(the number of rows is much larger than the number of columns),
$\boldsymbol{d}(t) \in \mathbb{C}^{K\times 1}$ is a vector
containing unknown amplitudes, and $1 \leq t \leq N$. A practically
important example of $\boldsymbol{A}$ is given in
Section~\ref{sec:App}. Since $\boldsymbol{A}$ is a tall matrix, the
covariance matrix of the signal is a low-rank matrix. Therefore,
such a signal is called low-rank. Matrix $\boldsymbol{A}$ can be
fully known, known up to a number of unknown parameters, or
completely unknown. In this , we study the second case where
matrix $\boldsymbol{A}$ has a known structure, but it contains $P$
number of unknown parameters $\boldsymbol{\Omega} \triangleq \left[
\omega_1, \cdots, \omega_P  \right]^T \in \mathbb{R}^{P\times 1}$
where $(\cdot)^T$ stands for the transposition operator.

Let the vector of the measurements $\boldsymbol{y} \in
\mathbb{C}^{N_y\times 1}$ be given by
\begin{eqnarray}
\boldsymbol{y}(t) \hspace{-2mm} &=& \hspace{-2mm} \boldsymbol{\Phi}
\left( \boldsymbol{x}(t) + \boldsymbol{w}(t) \right) \nonumber \\
&=& \hspace{-2mm} \boldsymbol{\Phi} \boldsymbol{x}(t) +
\boldsymbol{n}(t) \label{eq:measurement}
\end{eqnarray}
where $\boldsymbol{\Phi} \in \mathbb{R}^{N_y\times N_x}$ is the
measurement matrix with $N_y \leq N_x$. The additive noise
$\boldsymbol{w}(t) \in\mathbb{C}^{N_x\times 1}$ is assumed to have
the circularly-symmetric complex jointly-Gaussian distribution
$\mathcal{N}_C(\boldsymbol{0},\sigma^2 \boldsymbol{I}_{N_x})$ where
$\boldsymbol{I}_{N_x}$ is the identity matrix of size $N_x$ and
$\sigma^2$ is the noise power. No specific structure for the 
measurement matrix $\boldsymbol{\Phi}$ needs to be considered in 
our derivations. %For example, the elements of $\boldsymbol{\Phi}$ can
%be drawn independently from Gaussian distribution $\mathcal{N}
%(0,1/{N_y})$. 
It is because $\boldsymbol{\Phi}$ is assumed to be known at the 
signal reconstruction stage, and therefore, it is treated as a 
deterministic matrix in our derivations. As a result, irrespective 
to how $\boldsymbol{\Phi}$ is generated, the measurement noise 
$\boldsymbol{n}(t) \in \mathbb{C}^{N_y\times 1}$ has Gaussian 
distribution $\mathcal{N}_C (\boldsymbol{0}, \boldsymbol{R} )$ 
where $\boldsymbol{R} ~ \hspace{-1mm} = ~ \hspace{-1mm} \sigma^2 
\boldsymbol{\Phi} \boldsymbol{\Phi}^T$.

\section{Derivation of the CRB}
In this section, we derive the CRB for the signal model given by
\eqref{eq:matrix_sig_mod} and \eqref{eq:measurement}.

First, let the vector of parameters be defined as
\begin{equation} \label{UnknownParameters}
\boldsymbol{\vartheta} \triangleq \left[\bar{\boldsymbol{d}}^T(1),
\tilde{\boldsymbol{d}}^T(1), \cdots, \bar{\boldsymbol{d}}^T(N),
\tilde{\boldsymbol{d}}^T(N), \boldsymbol{\Omega}^T\right]^T
\end{equation}
where $\bar{\boldsymbol{d}}(t)$ and
$\tilde{\boldsymbol{d}}(t)$ represent the real and imaginary parts
of $\boldsymbol{d}(t)$, respectively.

The likelihood function of the compressed measurements
\eqref{eq:measurement} is given by
\begin{eqnarray} \label{Likelyhood}
&& \hspace{-10mm} p \left(\boldsymbol{y}(1), \cdots,
\boldsymbol{y}(N) \,|\, \boldsymbol{\vartheta} \right) =
\frac{1}{\pi^{N_y N} \left|
\boldsymbol{R} \right|^N} \nonumber \\
&& \hspace{-8mm} \times \exp\bigg\{- \sum_{t=1}^N
\left(\boldsymbol{y}(t)-\boldsymbol{B}\boldsymbol{d}(t)\right)^H
\boldsymbol{R}^{-1}
\left(\boldsymbol{y}(t)-\boldsymbol{B}\boldsymbol{d}(t)\right)\bigg\}
\end{eqnarray}
where $\boldsymbol{B} \triangleq \boldsymbol{\Phi} \boldsymbol{A}$ and
$(\cdot)^H$ stands for the Hermitian transposition operator. The
log-likelihood function can be found by taking the natural logarithm
of \eqref{Likelyhood} as
\begin{eqnarray}
LL \hspace{-2mm} & \triangleq & \hspace{-2mm} \ln p
\left(\boldsymbol{y}(1), \cdots, \boldsymbol{y}(N) \,|\,
\boldsymbol{\vartheta} \right) \nonumber \\
& = & \hspace{-2mm} - N_y N \ln\pi - N
\ln \left| \boldsymbol{R} \right| \nonumber \\
&& \hspace{-2mm} - \sum_{t=1}^N
\left(\boldsymbol{y}(t)-\boldsymbol{B}\boldsymbol{d}(t)\right)^H
\boldsymbol{R}^{-1}
\left(\boldsymbol{y}(t)-\boldsymbol{B}\boldsymbol{d}(t)\right).
\label{eq:LL}
\end{eqnarray}
For brevity, the notation $LL$ will be used in the rest of the 
to refer to the log-likelihood function \eqref{eq:LL}. The Fisher
information matrix (FIM) is given by
\begin{equation}
\boldsymbol{I}(\boldsymbol{\vartheta}) = E\left\{ \boldsymbol{\psi}
\boldsymbol{\psi}^T \right\}
\end{equation}
where $\boldsymbol{\psi} \triangleq \partial LL/\partial
\boldsymbol{\vartheta}$. The CRB covariance matrix for the vector of
parameters $\boldsymbol{\vartheta}$ is then given by
\begin{equation}
\text{CRB} \left( \boldsymbol{\vartheta} \right) =
\boldsymbol{I}^{-1}(\boldsymbol{\vartheta}).
\end{equation}

The derivatives of the $LL$ with respect to $\bar{\boldsymbol{d}}(t)$
and $\tilde{\boldsymbol{d}}(t)$ are given by
\begin{eqnarray}
\frac{\partial LL}{\partial\bar{\boldsymbol{d}}(t)} \hspace{-2mm}
&=& \hspace{-2mm} \boldsymbol{B}^H \boldsymbol{R}^{-1}
\boldsymbol{n}(t) + \left(\boldsymbol{n}^H(t) \boldsymbol{R}^{-1}
\boldsymbol{B}\right)^T
\nonumber \\
&=& \hspace{-2mm} 2 Re\left\{\boldsymbol{B}^H \boldsymbol{R}^{-1}
\boldsymbol{n}(t)\right\}%\label{eq:der_sbar}
\end{eqnarray}
and
\begin{eqnarray}
\frac{\partial LL}{\partial\tilde{\boldsymbol{d}}(t)} \hspace{-2mm}
&=& \hspace{-2mm} -j\boldsymbol{B}^H \boldsymbol{R}^{-1}
\boldsymbol{n}(t) + j\left(\boldsymbol{n}^H(t)
\boldsymbol{R}^{-1} \boldsymbol{B}\right)^T \nonumber \\
&=& \hspace{-2mm} 2 Im\left\{\boldsymbol{B}^H \boldsymbol{R}^{-1}
\boldsymbol{n}(t)\right\}%\label{eq:der_stilde}
\end{eqnarray}
where $Re\{\cdot\}$ and $Im\{\cdot\}$ stand for the real part and
imaginary part operators, respectively. Recall that
$\boldsymbol{n}(t)=\boldsymbol{y}(t)-\boldsymbol{B}\boldsymbol{d}(t)$
is the measurement noise introduced in \eqref{eq:measurement}.

Note that $\boldsymbol{A}$ has a known structure and contains $P$
unknown parameters $\omega_1, \cdots, \omega_P$. Therefore, the
derivative of the LL with respect to $\omega_p$ for $1\leq p\leq
P$ can be found as
\begin{eqnarray}
\frac{\partial LL}{\partial\omega_p} \hspace{-2mm} &=& \hspace{-2mm}
\sum_{t=1}^N \boldsymbol{d}^H(t)
\frac{\partial\boldsymbol{B}^H}{\partial\omega_p}
\boldsymbol{R}^{-1} \boldsymbol{n}(t) + \boldsymbol{n}^H(t)
\boldsymbol{R}^{-1} \frac{\partial\boldsymbol{B}}{\partial\omega_p}
\boldsymbol{d}(t) \nonumber \\
&=& \hspace{-2mm} 2 \sum_{t=1}^N Re\left\{\boldsymbol{d}^H(t)
\frac{\partial\boldsymbol{B}^H}{\partial\omega_p} \boldsymbol{R}^{-1}
\boldsymbol{n}(t) \right\} \nonumber \\
&=& \hspace{-2mm} 2 \sum_{t=1}^N Re\left\{\boldsymbol{d}^H(t)
\frac{\partial\boldsymbol{A}^H}{\partial\omega_p}
\boldsymbol{\Phi}^T \boldsymbol{R}^{-1} \boldsymbol{n}(t) \right\}.
\end{eqnarray}
The derivatives of the LL with respect to the whole vector
$\boldsymbol{\Omega}$ can be then written in matrix form as
\begin{equation}
\frac{\partial LL}{\partial\boldsymbol{\Omega}} = 2 \sum_{t=1}^N
Re\left\{\boldsymbol{D}^H(t) \boldsymbol{\Phi}^T \boldsymbol{R}^{-1}
\boldsymbol{n}(t)\right\}
\end{equation}
where the matrix $\boldsymbol{D}(t) \in \mathbb{C}^{N_x\times P}$ is
given by
\begin{eqnarray}
\boldsymbol{D}(t) \hspace{-2mm} &\triangleq& \hspace{-2mm} \left[
\frac{\partial\boldsymbol{A}}{\partial\omega_1} \boldsymbol{d}(t),
\cdots, \frac{\partial\boldsymbol{A}}{\partial\omega_P}
\boldsymbol{d}(t) \right] \nonumber \\
&=& \hspace{-2mm} \left[
\frac{\partial\boldsymbol{A}}{\partial\omega_1}, \cdots,
\frac{\partial\boldsymbol{A}}{\partial\omega_P} \right] \left(
\boldsymbol{I}_P \otimes \boldsymbol{d}(t) \right) \label{matrD}
\end{eqnarray}
with $\otimes$ standing for the Kronecker product.

To proceed, we use the following identities \cite{Stoica89:ML_CRB}.
For two arbitrary complex vectors $\boldsymbol{p}$ and
$\boldsymbol{q}$, we have
\begin{eqnarray}
Re(\boldsymbol{p})Re\left(\boldsymbol{q}^T\right) \hspace{-2mm}
&=& \hspace{-2mm}
\frac{1}{2}\left(Re\left(\boldsymbol{p}\boldsymbol{q}^T\right)+
Re\left(\boldsymbol{p}\boldsymbol{q}^H\right)\right) 
\label{eq:lemma_RpRq} \\
Im(\boldsymbol{p})Im\left(\boldsymbol{q}^T\right) 
\hspace{-2mm} &=& \hspace{-2mm}
-\frac{1}{2}\left(Re\left(\boldsymbol{p}\boldsymbol{q}^T\right)-
Re\left(\boldsymbol{p}\boldsymbol{q}^H\right)\right) 
\label{eq:lemma_IpIq} \\
Re(\boldsymbol{p})Im\left(\boldsymbol{q}^T\right) 
\hspace{-2mm} &=& \hspace{-2mm}
\frac{1}{2}\left(Im\left(\boldsymbol{p}\boldsymbol{q}^T\right)-
Im\left(\boldsymbol{p}\boldsymbol{q}^H\right)\right).
\label{eq:lemma_RpIq}
\end{eqnarray}
Using \eqref{eq:lemma_RpRq}, \eqref{eq:lemma_IpIq},
\eqref{eq:lemma_RpIq}, and the fact that for $1 \leq r,s \leq N$
\begin{eqnarray}
E \left\{ \boldsymbol{n}(r) \boldsymbol{n}^T(s) \right\}
\hspace{-2mm} &=&
\hspace{-2mm} \boldsymbol{0} \\
E \left\{ \boldsymbol{n}(r) \boldsymbol{n}^H(s) \right\}
\hspace{-2mm}
 &=& \hspace{-2mm} \delta_{r,s} \boldsymbol{R}
\end{eqnarray}
where $\delta_{r,s}$ denotes the Kronecker delta, we can compute
the submatrices of $\boldsymbol{I}(\boldsymbol{\vartheta})$ as
\begin{eqnarray}
E\left(\frac{\partial
LL}{\partial\bar{\boldsymbol{d}}(r)}\right)\left(\frac{\partial
LL}{\partial\bar{\boldsymbol{d}}(s)}\right)^T \hspace{-2mm} &=&
\hspace{-2mm} 2 Re\Big\{E \Big\{ \boldsymbol{B}^H \boldsymbol{R}^{-1} 
\boldsymbol{n}(r) \boldsymbol{n}^H(s) \boldsymbol{R}^{-1} 
\boldsymbol{B}\Big\}\Big\}\nonumber\\
\hspace{-2mm} &=& \hspace{-2mm} 2 Re\left\{ \boldsymbol{B}^H 
\boldsymbol{R}^{-1} \boldsymbol{B}\right\} \delta_{r,s} \\
E\left(\frac{\partial
LL}{\partial\bar{\boldsymbol{d}}(r)}\right)\left(\frac{\partial
LL}{\partial\tilde{\boldsymbol{d}}(s)}\right)^T \hspace{-2mm} &=&
\hspace{-2mm} -2 Im\Big\{E \Big\{ \boldsymbol{B}^H \boldsymbol{R}^{-1} 
\boldsymbol{n}(r) \boldsymbol{n}^H(s) \boldsymbol{R}^{-1} 
\boldsymbol{B}\Big\}\Big\}\nonumber\\
&=& \hspace{-2mm} - 2 Im\left\{ \boldsymbol{B}^H \boldsymbol{R}^{-1} 
\boldsymbol{B}\right\} \delta_{r,s} \\
E\left(\frac{\partial
LL}{\partial\bar{\boldsymbol{d}}(r)}\right)\left(\frac{\partial
LL}{\partial\boldsymbol{\Omega}}\right)^T \hspace{-2mm} &=&
\hspace{-2mm} 2 Re\Big\{E \Big\{ \boldsymbol{B}^H \boldsymbol{R}^{-1} 
\boldsymbol{n}(r) \sum_{t=1}^N \boldsymbol{n}^H(t)
\boldsymbol{R}^{-1} \boldsymbol{\Phi}
\boldsymbol{D}(t) \Big\}\Big\}\nonumber\\
&=& \hspace{-2mm} 2 Re\left\{\boldsymbol{B}^H \boldsymbol{R}^{-1} 
\boldsymbol{\Phi} \boldsymbol{D}(r)\right\}\\
E\left(\frac{\partial
LL}{\partial\tilde{\boldsymbol{d}}(r)}\right)\left(\frac{\partial
LL}{\partial\tilde{\boldsymbol{d}}(s)}\right)^T \hspace{-2mm} &=&
\hspace{-2mm} 2 Re\Big\{E \Big\{ \boldsymbol{B}^H \boldsymbol{R}^{-1} 
\boldsymbol{n}(r) \boldsymbol{n}^H(s) \boldsymbol{R}^{-1} 
\boldsymbol{B}\Big\}\Big\}\nonumber\\
&=& \hspace{-2mm} 2 Re\left\{\boldsymbol{B}^H \boldsymbol{R}^{-1} 
\boldsymbol{B}\right\} \delta_{r,s} \\
E\left(\frac{\partial
LL}{\partial\tilde{\boldsymbol{d}}(r)}\right)\left(\frac{\partial
LL}{\partial\boldsymbol{\Omega}}\right)^T \hspace{-2mm} &=&
\hspace{-2mm} - 2 Im\Big\{E \Big\{ \sum_{t=1}^N \boldsymbol{D}^H(t)
\boldsymbol{\Phi}^T \boldsymbol{R}^{-1} \boldsymbol{n}(t) 
\boldsymbol{n}^H(r) \boldsymbol{R}^{-1} \boldsymbol{B} \Big\}\Big\}^T 
\nonumber\\
&=& \hspace{-2mm} - 2 Im\Big\{ \boldsymbol{D}^H(r) \boldsymbol{\Phi}^T 
\boldsymbol{R}^{-1} \boldsymbol{B} \Big\}^T \nonumber\\
&=& \hspace{-2mm} 2 Im\left\{\boldsymbol{B}^H \boldsymbol{R}^{-1} 
\boldsymbol{\Phi}\boldsymbol{D}(r)\right\}\\
E\left(\frac{\partial
LL}{\partial\boldsymbol{\Omega}}\right)\left(\frac{\partial
LL}{\partial\boldsymbol{\Omega}}\right)^T \hspace{-2mm} &=&
\hspace{-2mm} 2 \sum_{t=1}^N \sum_{r=1}^N
Re \Big\{ E \Big\{ \boldsymbol{D}^H(t) \boldsymbol{\Phi}^T 
\boldsymbol{R}^{-1} \boldsymbol{n}(t)
\boldsymbol{n}^H(r) \boldsymbol{R}^{-1}
\boldsymbol{\Phi} \boldsymbol{D}(r) \Big\} \Big\} \nonumber\\
&=& 2 \sum_{t=1}^N Re\left\{\boldsymbol{D}^H(t)
\boldsymbol{\Phi}^T \boldsymbol{R}^{-1} \boldsymbol{\Phi}
\boldsymbol{D}(t)\right\}.
\end{eqnarray}
Then, $\boldsymbol{I}(\boldsymbol{\vartheta})$ can be found as
\begin{equation}
\boldsymbol{I}(\boldsymbol{\vartheta})=\left[
\begin{array}{cccccc}
\bar{\boldsymbol{H}} & -\tilde{\boldsymbol{H}} & & & &
\bar{\boldsymbol{\Delta}}_1\\
\tilde{\boldsymbol{H}} & \bar{\boldsymbol{H}} & & \boldsymbol{0}
& & \tilde{\boldsymbol{\Delta}}_1\\
& & \ddots & & & \vdots \\
\boldsymbol{0} & & & \bar{\boldsymbol{H}} & -\tilde{\boldsymbol{H}}
& \bar{\boldsymbol{\Delta}}_N\\
& & & \tilde{\boldsymbol{H}} & \bar{\boldsymbol{H}} &
\tilde{\boldsymbol{\Delta}}_N\\
\bar{\boldsymbol{\Delta}}_1^T & \tilde{\boldsymbol{\Delta}}_1^T &
\cdots & \bar{\boldsymbol{\Delta}}_N^T &
\tilde{\boldsymbol{\Delta}}_N^T & \boldsymbol{\Gamma}
\end{array} \right] \label{eq:I_theta_general}
\end{equation}
where $\bar{(\cdot)}$ and $\tilde{(\cdot)}$ stand for the real and
imaginary parts of a matrix, and
\begin{eqnarray}
\boldsymbol{H} \hspace{-2mm}  & \triangleq & \hspace{-2mm} 2 
\boldsymbol{B}^H \boldsymbol{R}^{-1} \boldsymbol{B} 
\label{eq:H_def} \\
\boldsymbol{\Delta}_r \hspace{-2mm}  & \triangleq & \hspace{-2mm} 
2 \boldsymbol{B}^H \boldsymbol{R}^{-1} \boldsymbol{\Phi}
\boldsymbol{D}(r) \label{eq:Delta_def} \\
\boldsymbol{\Gamma} \hspace{-2mm}  & \triangleq & \hspace{-2mm} 2
\sum_{t=1}^N Re\left\{\boldsymbol{D}^H(t) \boldsymbol{\Phi}^T
\boldsymbol{R}^{-1} \boldsymbol{\Phi} \boldsymbol{D}(t) \right\}.
\label{eq:Gamma_def}
\end{eqnarray}

It is shown in \cite{Stoica89:ML_CRB} that for FIM with the structure
given in \eqref{eq:I_theta_general}, the CRB covariance matrix for
$\boldsymbol{\Omega}$ is given by
\begin{equation}
\text{CRB}^{-1} \left( \boldsymbol{\Omega} \right) =
\boldsymbol{\Gamma} - \sum_{t=1}^N Re \left\{
\boldsymbol{\Delta}_t^H \boldsymbol{H}^{-1} \boldsymbol{\Delta}_t
\right\}. \label{eq:CRB_inv_Omega}
\end{equation}
Using \eqref{eq:H_def}--\eqref{eq:CRB_inv_Omega}, the CRB for
$\boldsymbol{\Omega}$ can be found in closed-form as
\begin{eqnarray}
\hspace{-6mm} \text{CRB}^{-1} \left( \boldsymbol{\Omega} \right)
\hspace{-2mm} &=& \hspace{-2mm} 2 \sum_{t=1}^N Re \Big\{
\boldsymbol{D}^H(t) \boldsymbol{\Phi}^T
\boldsymbol{R}^{-1} \Big( \boldsymbol{I}_{N_y} \nonumber \\
&& \hspace{-2mm} - \boldsymbol{B} \left( \boldsymbol{B}^H
\boldsymbol{R}^{-1} \boldsymbol{B} \right)^{-1} \boldsymbol{B}^H
\boldsymbol{R}^{-1} \Big) \boldsymbol{\Phi}\boldsymbol{D}(t)
\Big\}. \label{eq:CRB_inv_Omega_smpl}
\end{eqnarray}

Given $\boldsymbol{I}^{-1}(\boldsymbol{\vartheta})$, the 
covariance matrix of any unbiased estimator of $\boldsymbol{x}
(t)$, i.e., $\boldsymbol{C}_{\hat{\boldsymbol{x}}(t)}$, satisfies 
the inequality
\begin{equation}
\boldsymbol{C}_{\hat{\boldsymbol{x}}(t)} - \frac{\partial
\boldsymbol{x}(t)} {\partial \boldsymbol{\vartheta}}
\boldsymbol{I}^{-1} (\boldsymbol{\vartheta}) \frac{\partial
\boldsymbol{x}^H(t)} {\partial\boldsymbol{\vartheta}} \geq
\boldsymbol{0}. \label{eq:CRB_ineq}
\end{equation}
The signal $\boldsymbol{x}(t)$ can be written as
\begin{equation}
\boldsymbol{x}(t) = \boldsymbol{A} \boldsymbol{d}(t) =
\boldsymbol{A} \bar{\boldsymbol{d}}(t) + j \boldsymbol{A}
\tilde{\boldsymbol{d}}(t).
\end{equation}
The derivative of $\boldsymbol{x}(t)$ with respect to the vector 
of unknown parameters $\boldsymbol{\vartheta}$ is given by
\begin{equation}
\frac{\partial\boldsymbol{x}(t)}{\partial\boldsymbol{\vartheta}} =
\left[ \boldsymbol{e}_t \otimes \left[ \boldsymbol{A},
j\boldsymbol{A} \right], \boldsymbol{D}(t) \right]
\label{eq:CRB_def1}
\end{equation}
where $\boldsymbol{e}_t$ is a row vector of length $N$ with all its
elements equal to zero except for the $t$-th element which is equal
to $1$. Finally, by summing over the diagonal elements of
\eqref{eq:CRB_ineq}, we obtain
\begin{equation}
E \left\{ \| \hat{\boldsymbol{x}}(t) - \boldsymbol{x}(t) \|^2
\right\} \geq \text{Tr} \left\{ \frac{\partial\boldsymbol{x}(t)}
{\partial\boldsymbol{\vartheta}} \boldsymbol{I}^{-1}
(\boldsymbol{\vartheta}) \frac{\partial\boldsymbol{x}^H(t)}
{\partial\boldsymbol{\vartheta}} \right\}.\label{eq:CRB_def}
\end{equation}
Similar to \eqref{eq:CRB_inv_Omega_smpl}, the results in
\eqref{eq:CRB_def1} and \eqref{eq:CRB_def} can also be regarded as
closed-form as they can be used for analysis, fast computations, 
and getting insights without requiring Monte-Carlo simulations as 
shown in the next section. It is worth noting that the derived CRB 
(especially \eqref{eq:CRB_inv_Omega_smpl}) can be also used for 
selecting/optimizing the measurement matrix $\boldsymbol{\Phi}$ as 
the CRB depends on a specific selection of $\boldsymbol{\Phi}$.

\section{Minimum Number of Compressed Samples}
\label{sec:min_Ny}
In this section, we show that if the number of compressed samples is
less than or equal to the number of sources ($N_y \leq K$), the
FIM $\boldsymbol{I}(\boldsymbol{\vartheta})$ is singular. It is
shown in \cite{Stoica01:Singular_FIM} that a singular FIM means that
unbiased estimation of the entire parameter vector with finite
variance is impossible.

Let us start with the case that $N_y < K$. In this case, we have
$rank \left\{ \boldsymbol{B} \right\} < K$ since $\boldsymbol{B} \in
\mathbb{C}^{N_y \times K}$. As a result, we also have $rank \left\{
\boldsymbol{H} \right\} < K$ (see \eqref{eq:H_def}), and therefore,
$\boldsymbol{H}$ is singular. Consequently, there exists a nonzero
vector $\boldsymbol{u} = \bar{\boldsymbol{u}} + j
\tilde{\boldsymbol{u}} \in \mathbb{C}^{K \times 1}$ such that
$\boldsymbol{H}\boldsymbol{u} = \boldsymbol{0}$. Therefore, $\left(
\bar{\boldsymbol{H}} + j \tilde{\boldsymbol{H}} \right) \left(
\bar{\boldsymbol{u}} + j \tilde{\boldsymbol{u}} \right) =
\boldsymbol{0}$, which can be written in matrix form as
\begin{equation}
\left[
\begin{array}{cc}
\bar{\boldsymbol{H}} & -\tilde{\boldsymbol{H}} \\
\tilde{\boldsymbol{H}} & \bar{\boldsymbol{H}}
\end{array} \right]
\left[
\begin{array}{c}
\bar{\boldsymbol{u}} \\
\tilde{\boldsymbol{u}}
\end{array} \right] = \boldsymbol{0}.
\end{equation}
Let $\boldsymbol{v} \triangleq \left[ \bar{\boldsymbol{u}}^T,
\tilde{\boldsymbol{u}}^T, \boldsymbol{0} \right]^T \in
\mathbb{R}^{(2NK+P) \times 1}$. Finally, using
\eqref{eq:I_theta_general}, we have $\boldsymbol{v}^T
\boldsymbol{I}(\boldsymbol{\vartheta}) \boldsymbol{v} =
\boldsymbol{0}$, which means that
$\boldsymbol{I}(\boldsymbol{\vartheta})$ has a zero eigenvalue,
and therefore, it is singular. For the case that $N_y = K$, if
$rank \left\{ \boldsymbol{B} \right\} < K$, the singularity of
the FIM follows from the discussion above.

Now, consider the case that
$\boldsymbol{B}$ is full-rank. Thus, $\boldsymbol{H}$ is invertible.
Consider the structure of $\boldsymbol{I}(\boldsymbol{\vartheta})$
in \eqref{eq:I_theta_general} and let the block of all the real and
imaginary parts of $\boldsymbol{H}$ be denoted by $\boldsymbol{T}$.
It is shown in \cite{Stoica89:ML_CRB} that for an invertible matrix
$\boldsymbol{H}$, matrix $\boldsymbol{T}$ is also invertible. The
Schur complement of $\boldsymbol{T}$ denoted by
$\boldsymbol{I}(\boldsymbol{\vartheta})/\boldsymbol{T}$ is equal to
the inverse of the CRB covariance matrix for $\boldsymbol{\Omega}$
as given in \eqref{eq:CRB_inv_Omega_smpl}. Matrix $\boldsymbol{B}$
is invertible since it is square and full-rank. Therefore, we have
\begin{equation}
\boldsymbol{I}_{N_y} - \boldsymbol{B} \left( \boldsymbol{B}^H
\boldsymbol{R}^{-1} \boldsymbol{B} \right)^{-1} \boldsymbol{B}^H
\boldsymbol{R}^{-1} = \boldsymbol{0}.
\end{equation}
As a result, $\boldsymbol{I}(\boldsymbol{\vartheta})/\boldsymbol{T}
= \boldsymbol{0}$ (see \eqref{eq:CRB_inv_Omega_smpl}). According to
the rank additivity formula \cite{Zhang05:Schur}, we have
\begin{equation}
rank \left\{ \boldsymbol{I}(\boldsymbol{\vartheta}) \right\} = rank
\left\{ \boldsymbol{T} \right\} + rank \left\{
\boldsymbol{I}(\boldsymbol{\vartheta})/\boldsymbol{T} \right\} =
rank \left\{ \boldsymbol{T} \right\}.
\end{equation}
Therefore, $\boldsymbol{I}(\boldsymbol{\vartheta})$ is
rank-deficient or equivalently singular.

\theoremstyle{remark}
\newtheorem*{rem_min_CS_num}{Remark}
\begin{rem_min_CS_num}
As shown above, if the number of compressed samples is less than or
equal to the number of sources, the FIM is necessarily singular. However,
if the number of compressed samples increases, it does not necessarily
result in a non-singular FIM for a few more samples. Thus, the
converse does not hold in general. The minimum number of compressed
samples for satisfactory performance depends on a specific performance
criterion and the estimation method used. For example, the minimum
number of compressed samples can be chosen to bound the probability
of a subspace swap \cite{Scharf_13:CS_sub_swap} or to bound the
error of signal subspace estimation \cite{Romberg14:CS_continuum}.
The required number of compressed samples can also be studied from a
geometric point of view \cite{Wakin09:Dim_reduction}.
\end{rem_min_CS_num}

\section{Application Examples}
\label{sec:App}
For the problems of DOA and spectral estimation, $\boldsymbol{d}(t)$
consists of the amplitudes of $K$ number of sources at time instant
$t$. The number of parameters in $\boldsymbol{A}$ is also equal to
the number of sources, i.e., $P = K$. Furthermore, $\boldsymbol{A}$
has the structure given by
\begin{equation}
\boldsymbol{A} \triangleq \left[ \boldsymbol{a}\left( \omega_1 \right),
\cdots, \boldsymbol{a}\left( \omega_K \right) \right]
\end{equation}
where $\boldsymbol{a}\left( \omega_k \right)$ for $1 \leq k \leq K$
is the steering vector corresponding to the $k$-th source. Let us
define $\boldsymbol{c}(\omega)$ as the derivative of $\boldsymbol{a}
(\omega)$ with respect to $\omega$, i.e., $\boldsymbol{c} (\omega) 
\triangleq d\boldsymbol{a}(\omega) / d\omega$. Then, $\boldsymbol{D}
(t)$ given by \eqref{matrD} can be simplified to
\begin{eqnarray}
\boldsymbol{D}(t) \hspace{-2mm} &=& \hspace{-2mm} \left[
\boldsymbol{c}\left( \omega_1 \right)d_1(t), \cdots,
\boldsymbol{c}\left( \omega_K \right)d_K(t) \right] \nonumber \\
&=& \hspace{-2mm} \left[ \boldsymbol{c}\left( \omega_1 \right),
\cdots, \boldsymbol{c}\left( \omega_K \right) \right]
\text{diag}\left\{ \boldsymbol{d}(t) \right\}
\end{eqnarray}
where $d_k(t)$ is the $k$-th element of $\boldsymbol{d}(t)$ and the
$\text{diag}\left\{ \cdot \right\}$ operator converts a vector into
a diagonal matrix.

\section{Numerical Results}
In this section, the application of the derived CRB formulas for the
problem of DOA estimation is illustrated. Our goal is to investigate
the performance bounds for unbiased estimators when the signal is
compressed at different rates.

Consider $K = 11$ equally spaced sources impinging on a uniform
linear array of $N_x = 50$ antenna elements from directions
$\omega_1 = 20\,^{\circ} \times (\pi/180), \omega_2 = 23\,^{\circ}
\times (\pi/180), \cdots, \omega_{11} = 50\,^{\circ} \times
(\pi/180)$. The steering vector of the array
$\boldsymbol{a}(\omega)$ is given by
\begin{equation}
\boldsymbol{a}(\omega) \triangleq \left[1, e^{-j2 \pi (d/\lambda)
\sin(\omega)}, \cdots, e^{-j2 \pi (N_x-1)(d/\lambda)
\sin(\omega)}\right]^T
\end{equation}
where $d$ is the interelement spacing of the array and $\lambda$ is
the wavelength of the plane wave impinging on the array. In our
numerical example, $d/\lambda$ is set to $0.5$. The number of
snapshots is also set to $N = 10$. Each source vector
$\boldsymbol{d}(t)$ is considered to be independent from the source
vectors at other time instances and is drawn from the
circularly-symmetric complex jointly-Gaussian distribution
$\mathcal{N}_C(\boldsymbol{0},\sigma_s^2 \boldsymbol{I}_K)$. The
signal-to-noise ratio (SNR) is set to $\text{SNR} \triangleq
10\log_{10} \left( \sigma_s^2 / \sigma^2 \right) = 10$~dB. The
source vectors are drawn once and kept unchanged.

\begin{figure}[t]
\psfrag{CRB}{CRB (dB)} \psfrag{Ny}{Number of compressed samples
($N_y$)} \psfrag{a}{\scriptsize{$N_y = 12$}} \hspace{0mm}
    \includegraphics[width=28cm]{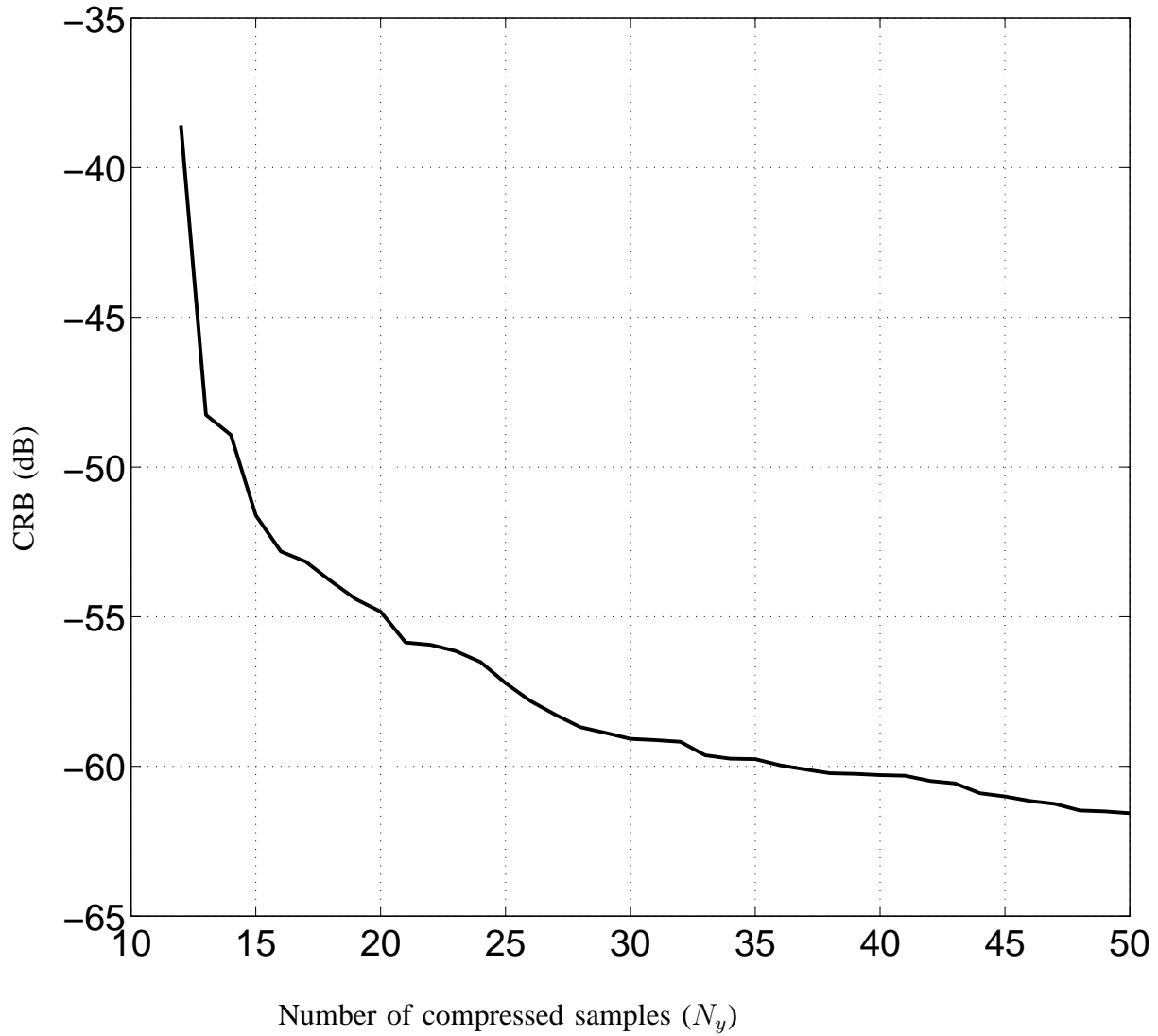}
\caption{CRB for estimating $\omega_6 = 35\,^{\circ} \times
(\pi/180)$. \label{fig:CRB_Ny}} 
\end{figure}

Fig.~\ref{fig:CRB_Ny} shows the CRB for estimating $\omega_6 =
35\,^{\circ} \times (\pi/180)$ versus the number of compressed
samples $N_y$. For the case when $N_y = N_x = 50$, the measurement
matrix $\boldsymbol{\Phi}$ is set to the identity matrix. Then,
$\boldsymbol{\Phi}$ is initialized for $N_y = 49$ by drawing samples
from the Gaussian distribution $\mathcal{N}(0,1/{49})$. For the rest
of $N_y$ values, the first $N_y$ rows of the initial matrix
$\boldsymbol{\Phi}$ are scaled by $\sqrt{49/N_y}$ and used to obtain
the CRB.

As expected, it can be seen in Fig.~\ref{fig:CRB_Ny} that the CRB
increases as the number of compressed samples $N_y$ reduces. The
minimum number of compressed samples is set to $N_y = 12$ which is
equal to the number of sources plus one ($K+1$). As shown in
Section~\ref{sec:min_Ny}, if the number of compressed samples is
equal to or less than the number of sources, there can be no 
unbiased estimator with a finite estimation variance. Otherwise, 
if the CRB exists, there also exist estimators 
\cite{MSH11:Nested_LS} that achieve it.

\newpage

\newpage
\section{Conclusion}
\label{sec:conclude} The class of signals fitting a low-rank 
signal model has been considered in this . Such signals are
inherently sparse according to the signal model and can be 
recovered from compressed 

\newpage \hspace{-0.5cm}
measurements. We have studied the
performance bounds for unbiased estimators of parameters of such 
low-rank signal model from compressed samples. The Cram\'{e}r-Rao
bound has been derived for a generic low-rank model and it has 
been shown that the number of compressed samples needs to be at
least larger than the number of sources for the existence of an 
unbiased estimator with finite variance. Furthermore, the 
applications to DOA and spectral estimation have been considered. 
Numerical examples have been also given to illustrate the effect 
of compression on the CRB. It has been shown how the CRB increases 
until the point where the number of compressed samples is larger 
than the number of sources. For lower number of compressed 
samples, the CRB becomes unbounded.

\bibliographystyle{IEEEtran}

\begin{comment}
\bibliography{references}

\begin{thebibliography}{10}
\providecommand{\url}[1]{#1} \csname url@rmstyle\endcsname
\providecommand{\newblock}{\relax} \providecommand{\bibinfo}[2]{#2}
\providecommand\BIBentrySTDinterwordspacing{\spaceskip=0pt\relax}
\providecommand\BIBentryALTinterwordstretchfactor{4}
\providecommand\BIBentryALTinterwordspacing{\spaceskip=\fontdimen2\font
plus \BIBentryALTinterwordstretchfactor\fontdimen3\font minus
  \fontdimen4\font\relax}
\providecommand\BIBforeignlanguage[2]{{%
\expandafter\ifx\csname l@#1\endcsname\relax
\typeout{** WARNING: IEEEtran.bst: No hyphenation pattern has been}%
\typeout{** loaded for the language `#1'. Using the pattern for}%
\typeout{** the default language instead.}%
\else \language=\csname l@#1\endcsname \fi #2}}

\bibitem{Donoho06:CS}
D.~Donoho, ``Compressed sensing,'' \emph{IEEE Trans. Inform. Theory}, 
vol.~52, no.~4, pp.~1289--1306, Apr.~2006.

\bibitem{Candes06:ExactCS}
E.~J.~Cand\`{e}s, J.~Romberg, and T.~Tao, ``Robust uncertainty principles: 
\textsc{E}xact signal reconstruction from highly incomplete frequency
information,'' \emph{IEEE Trans. Inform. Theory}, vol.~52, no.~2,
pp.~489--509, Feb.~2006.

\bibitem{Candes06:CS}
E.~J.~Cand\`{e}s and T.~Tao, ``Near optimal signal recovery from random
projections: \textsc{U}niversal encoding strategies?'' \emph{IEEE 
Trans. Inform. Theory}, vol.~52, no.~12, pp.~5406--5425, Dec.~2006.
	
\bibitem{Candes08:CS}
E.~J.~Cand\`{e}s and M.~B.~Wakin, ``An introduction to compressive 
sampling,'' \emph{IEEE Signal Process. Mag.}, vol.~25, no.~2, 
pp.~21-–30, Mar.~2008.

\bibitem{Laska07:CS}
J.~N.~Laska, S.~Kirolos, M.~F.~Duarte, T.~S.~Ragheb, R.~G.~Baraniuk,
and Y.~Massoud, ``Theory and implementation of an analog-to-information
converter using random demodulation,'' in \emph{Proc. IEEE ISCAS},
New Orleans, LA, May~2007, pp.~1959-–1962.

\bibitem{Vorobyov11:CS}
O.~Taheri and S.~A.~Vorobyov, ``Segmented compressed sampling for 
analog-to-information conversion: Method and performance analysis,'' 
\emph{IEEE Trans. Signal Process.}, vol.~59, no.~2, pp.~554--572, 
Feb.~2011.

\bibitem{Vorobyov14:CS}
H.~Fang, S.~A.~Vorobyov, and H.~Jiang, ``Performance limits of 
segmented compressive sampling: Correlated samples versus bits,'' 
\emph{preprint arXiv: 1411.5178}, 2014.

\bibitem{Xam1}
M.~Mishali and Y.~C.~Eldar, ``Sub-Nyquist sampling: Bridging theory 
and practice,'' \emph{IEEE Signal Process. Mag.}, vol.~28, no.~6, 
pp.~98--124, Nov.~2011.

\bibitem{Xam2}
M.~Mishali, Y.~C.~Eldar, and A.~Elron, ``Xampling: Signal acquisition 
and processing in union of subspaces,'' \emph{IEEE Trans. Signal 
Process.}, vol.~59, no.~10, pp.~4719--4734, Oct.~2011. 

\bibitem{FigNowak07}
M.~A.~T. Figueiredo, R.~D. Nowak, and S.~J. Wright, ``Gradient
projection for sparse reconstruction: Application to compressed 
sensing and other inverse problems,'' \emph{IEEE J. Select. Topics 
Signal Process.}, vol.~1, no.~4, pp.~586--597, Dec.~2007.

\bibitem{Al1}
E.~J.~Cand\`{e}s and T. Tao, ``Decoding by linear programming,'' 
\emph{IEEE Trans. Inf. Theory}, vol.~51, no.~12, pp.~4203-–4215, 
Dec.~2005.

\bibitem{Al2}
E.~J.~Cand\`{e}s J.~Romberg, and T.~Tao, ``Stable signal recovery 
from incomplete and inaccurate measurements,'' \emph{Commun. Pure 
Appl. Math.}, vol.~59, pp.~1207-–1223, Aug.~2006.

\bibitem{Al3}
J.~Haupt and R.~Nowak, ``Signal reconstruction from noisy random 
projections,'' \emph{IEEE Trans. Inf. Theory}, vol.~52, no.~9, 
pp.~4036-–4048, Sep.~2006.

\bibitem{Al4}
J.~A.~Tropp and A.~C.~Gilbert, ``Signal recovery from random 
measurements via orthogonal matching pursuit,'' \emph{IEEE Trans. 
Inf. Theory}, vol.~53, no.~12, pp.~4655-4666, Dec.~2007.

\bibitem{Al5}
D.~Needell and J.~A.~Tropp, ``CoSaMP: Iterative signal recovery 
from incomplete and inaccurate samples,'' \emph{Appl. and Comput. 
Harmonic Analysis}, vol.~26, no.~3, pp.~301--321, May~2009.

\bibitem{Al6}
E.~J.~Cand\`{e}s, M.~B.~Wakin, and S.~P.~Boyd, ``Enhancing 
sparsity by reweighted $l_1$ minimization,'' \emph{J. Fourier 
Anal. Appl.}, vol.~14, no.~5, pp.~877–-905, Dec.~2008.

\bibitem{Al7}
O.~Taheri and S.~A.~Vorobyov, ``Reweighted $l_1$-norm penalized 
LMS for sparse channel estimation and its analysis,'' \emph{Signal 
Process.}, vol.~104, pp.~70--79, May~2014.

\bibitem{MSH11:Nested_LS}
M.~Shaghaghi and S.~A. Vorobyov, ``Improved model-based spectral
compressive sensing via nested least squares,'' in \emph{Proc. IEEE 
Int. Conf. Acoustics, Speech, Signal Process. (ICASSP)}, Prague, 
Czech Republic, May~2011, pp.~3904--3907.

\bibitem{Scharf11:Basis_Mismatch}
Y.~Chi, L.~L. Scharf, A.~Pezeshki, and A.~R. Calderbank,
``Sensitivity to basis mismatch in compressed sensing,'' \emph{IEEE 
Trans. Signal Process.}, vol.~59, no.~5, pp.~2182--2195, May~2011.

\bibitem{VanTrees}
H.~L.~Van~Trees, K.~L.~Bell, and Z.~Tian, \emph{Detection, 
Estimation, and Modulation Theory.}\hskip
1em plus 0.5em minus 0.4em\relax USA: John Wiley \& Sons, 2005.

\bibitem{Eldar10:Sparse_CRB}
Z.~Ben-Haim and Y.~C. Eldar, ``The \textsc{C}ram\'{e}r-\textsc{R}ao
bound for estimating a sparse parameter vector,'' \emph{IEEE Trans. 
Signal Process.}, vol.~58, no.~6, pp.~3384--3389, Jun.~2010.

\bibitem{Gorman90:Constrained_CRB}
J.~D. Gorman and A.~O. Hero, ``Lower bounds for parametric
estimation with constraints,'' \emph{IEEE Trans. Inform. Theory}, 
vol.~26, no.~6, pp.~1285--1301, Nov.~1990.

\bibitem{Stoica98:Constrained_CRB}
P.~Stoica and B.~C. Ng, ``On the \textsc{C}ram\'{e}r-\textsc{R}ao
bound under parametric constraints,'' \emph{IEEE Signal Process. 
Lett.}, vol.~5, no.~7, pp.~177--179, Jul.~1998.

\bibitem{Nehorai11:LowRank_CRB}
G.~Tang and A.~Nehorai, ``Lower bounds on the mean-squared error of
low-rank
  matrix reconstruction,'' \emph{IEEE Trans. Signal Process.}, vol.~59, no.~10,
  pp.~4559--4571, Oct.~2011.

\bibitem{Scharf13:CS_Param_CRB}
P.~Pakrooh, L.~L. Scharf, A.~Pezeshki, and Y.~Chi, ``Analysis of
fisher
  information and the \textsc{C}ram\'{e}r-\textsc{R}ao bound for nonlinear
  parameter estimation after compressed sensing,'' in \emph{Proc. IEEE Int.
  Conf. Acoustics, Speech, Signal Process. (ICASSP)}, Vancouver, BC, May~2013, \mbox{pp.~6630--6634}.

\bibitem{Ramasamy12:CS_CRB_asilomar}
D.~Ramasamy, S.~Venkateswaran, and U.~Madhow, ``Compressive
estimation in
  \textsc{AWGN}: general observations and a case study,'' in \emph{Proc. 46th
  Asilomar Conf. Signals, Syst., Comput. (ASILOMAR)}, 2012, \mbox{pp.~953--957}.

\bibitem{Ramasamy14:CS_param_TSP}
D.~Ramasamy, S.~Venkateswaran, and U.~Madhow, ``Compressive parameter 
estimation in \textsc{AWGN},'' \emph{IEEE Trans. Signal Process.}, vol.~62, 
no.~8, pp.~2012--2027, Apr.~2014.

\bibitem{Stoica89:ML_CRB}
P.~Stoica and A.~Nehorai, ``\textsc{MUSIC}, maximum likelihood, and
  \textsc{C}ramer-\textsc{R}ao bound,'' \emph{IEEE Trans. Acoust., Speech,
  Signal Process.}, vol.~37, no.~5, pp.~720--741, May~1989.

\bibitem{Stoica01:Singular_FIM}
P.~Stoica and T.~L. Marzetta, ``Parameter estimation problems with
singular
  information matrices,'' \emph{IEEE Trans. Signal Process.}, vol.~49, no.~1,
  pp.~87--90, Jan.~2001.

\bibitem{Zhang05:Schur}
F.~Zhang, \emph{The Schur Complement and Its Applications.}\hskip
1em plus
  0.5em minus 0.4em\relax Berlin, Germany: Springer-Verlag, 2005.

\bibitem{Scharf_13:CS_sub_swap}
P.~Pakrooh, A.~Pezeshki, and L.~L. Scharf, ``Threshold effects in
parameter
  estimation from compressed data,'' in \emph{Proc. 1st IEEE Global Conference
  on Signal and Information Processing (GlobalSIP)}, Austin, TX, Dec.~2013, pp.~997--1000.

\bibitem{Romberg14:CS_continuum}
W.~Mantzel and J.~Romberg, ``Compressed subspace matching on the
continuum,''
  \emph{preprint arXiv:1407.5234}, 2014.

\bibitem{Wakin09:Dim_reduction}
R.~Baraniuk and M.~Wakin, ``Random projections of smooth
manifolds,''
  \emph{Found. Comput. Math.}, vol.~9, no.~1, pp.~51--77, Feb.~2009.

\end{thebibliography}
\end{comment}

\end{document}